\documentstyle[amsfonts,11pt]{article}

\newtheorem{teor}{Theorem}[section]

\newtheorem{remar}[teor]{Remark}

\newtheorem{lemma}[teor]{Lemma}

\newcommand{\fdim}{\hspace*{\fill}$\Box$}
\newcommand{\dimostr}{{\bf Proof: }}

\newcommand{\real}{\Bbb{R}}

\newcommand{\complex}{\Bbb{C}}

\newcommand{\natur}{\Bbb{N}}

\newcommand{\K}{K\"{a}hler}

\begin{document}

\noindent
\centerline {\LARGE\bf  Radial balanced metrics
on the unit disk}





\vspace{0.5cm}

\centerline{\small Antonio Greco and Andrea Loi}
\centerline{\small Dipartimento di Matematica e Informatica
-- Universit\`{a} di Cagliari}
\centerline{\small Via Ospedale 72, 09124 Cagliari -- Italy}
\centerline{\small e-mail : greco@unica.it, loi@unica.it}

\vskip 0.5cm

\vskip 0.3cm

\begin{abstract}
\noindent
Let   $\Phi$ be a  strictly plurisubharmonic and   radial function
on the unit disk ${\cal D}\subset {\complex}$
and let $g$ be the \K\ metric associated to the \K\ form
$\omega =\frac{i}{2}\partial\bar\partial\Phi$.
We prove that  if 
$g$ is $g_{eucl}$-balanced  of height $3$  (where $g_{eucl}$
is the standard  Euclidean metric  on ${\complex}={\real}^2$), 
and the function
$h(x)=e^{-\Phi (z)}$,  $x=|z|^2$, extends to  an entire analytic function on ${\real}$, 
then  $g$ equals the hyperbolic metric.
The proof of  our result
is based on a interesting  characterization of the  function 
$f(x)=1-x$.

\vskip 0.3cm


\noindent
{\it{Subj.Class}}: 53D05, 53C55, 58C25, 58F06.

\end{abstract}

\section{Introduction and statement of the main results}

\noindent

Let $\Phi :M\rightarrow {\real}$ be a strictly plurisubharmonic function
on a $n$-dimensional complex manifold $M$ and  
let $g_0$ be a \K\ metric on $M$.
Denote by ${\cal H}=L_{hol}^2(M, e^{-\Phi}\frac{\omega_0^n}{n!})$
the separable  complex Hilbert space
consisting of holomorphic  functions 
$\varphi$ on $M$ such that
\begin{equation}\label{piomega}
\langle \varphi, \varphi\rangle=
\int_M e^{-\Phi}|\varphi|^2\frac{\omega_0^n}{n!}<\infty,
\end{equation}
where $\omega_0$ is the \K\ form
associated to  the \K\ metric $g_0$
(this means that $\omega_0 (X, Y)=g_0(JX, Y)$, 
for all vector fields $X, Y$ on $M$, where 
$J$ is the complex structure of $M$).
Assume for each point $x\in M$
there exists $\varphi\in {\cal H}$
non-vanishing at $x$.
Then,  one can consider the
following holomorphic map into the
$N$-dimensional ($N\leq\infty$) complex projective space:
\begin{equation}\label{psiglob}
\varphi_{\Phi}:M\rightarrow {\complex}P^N:
x\mapsto [\varphi_0(x),\dots ,\varphi_N(x)],
\end{equation}
where $\varphi_j,\,\, j=0,\dots ,N,$
is a orthonormal basis for
${\cal H}$.
In the case $N=\infty$,
${\complex}P^{\infty}$
denote the quotient space of $l^2({\complex})$
(the space of sequences $z_j$ such that
$\sum_{j=1}^{\infty}|z_j|^2<\infty$),
where two sequences $z_j$ and $w_j$
are equivalent iff there exists $\lambda\in {\complex}^*={\complex}\setminus \{0\}$
such that $w_j=\lambda z_j , \forall j$.

Let $g$ be the \K\ metric associated to the \K\ form
$\omega =\frac{i}{2}\partial\bar\partial\Phi$ (and 
so $\Phi$ is a \K\ potential for $g$).
We say that the metric  $g$
is {\em $g_0$-balanced
of height $\alpha$}, $\alpha >0$, 
if $\varphi_{\Phi}^*g_{FS}=\alpha g$,
or equivalently
\begin{equation}\label{projind}
\varphi_{\Phi}^*\omega_{FS}=\alpha\omega,
\end{equation}
where $g_{FS}$ is the Fubini--Study metric on
${\complex}P^N$ and $\omega_{FS}$
its associated \K\ form,
namely 
$$\omega_{FS}=\frac{i}{2}
\partial\bar{\partial}\log\displaystyle
\sum _{j=0}^{N}|Z_{j}|^{2},$$
for  a  homogeneous coordinate
system
$[Z_0,\dots , Z_{N}]$ of 
${\complex}P^{N}$ (note that this definition
is independent from the choice of
the orthonormal
basis).
Therefore, if  $g$ is a $g_0$-balanced  
metric of height $\alpha$,
then  $\alpha g$ is projectively induced
via the map (\ref{psiglob}) (we refer the reader to the 
seminal paper \cite{ca} for more details on 
projectively induced metrics).
In the case
a metric $g$ is $g$-balanced, i.e. $g=g_0$,
one simply call $g$
a {\em balanced} metric.

The study of balanced and $g_0$-balanced
metrics is a very fruitful area of research both from
mathematical and physical point of view
(see \cite{arlcomm}, \cite{cgr3}, \cite{cgr4}, \cite{do3},  \cite{me}, \cite{mebal}, \cite{albounded},
\cite{aleps}, \cite{regscal}  and 
\cite{albergbal}).
The map $\varphi_{\Phi}$ was
introduced by J. Rawnsley \cite{ra}
in the context of quantization of \K\ manifolds
and it is often referred to as
the {\em coherent states map}.

Notice that one can easily  give an alternative
definition of  balanced metrics (not involving
projectively induced \K\ metrics)
in terms of the  reproducing kernel 
of the Hilbert space ${\cal H}$. Nevertheless
the  defintion given here is motivated by the recent results
on compact manifolds. In fact,   it can be easily extended  
 to the  case when $(M, \omega )$ is a polarized  compact
 \K\ manifold, with polarization $L$, 
i.e.,  $L$ is   a holomorphic line bundle 
$L$ over~$M$,
such that   $c_1(L)=[ \omegaÊ]$
(see e.g. \cite{arlquant} and \cite{arlfirst}
for details).
In the quantum mechanics terminology the bundle $L$
is called the {\em quantum line bundle} and the pair  $(L, h)$
a {\em  geometric quantization} of $(M, \omega)$.
The problem of the existence and uniqueness   of balanced metrics on a given \K\
class  of a compact manifold $M$ was solved by S. Donaldson  \cite{do}
when the group of  biholomorphisms of $M$ which lifts to the quantum line bundle $L$
modulo the ${\complex}^*$ action is finite and by C. Arezzo and the second author in the 
general case (see also \cite{ma}).

Nevertheless, many basic	 and important  questions on the existence and  uniqueness
of balanced metrics on noncompact manifolds are still open.
For example, it is unknown if there exists
a complete balanced metric on ${\complex}^n$
different from the euclidean metric.
The case of  $g_0$-balanced metric on
${\complex}^n$, where $g_0=g_{eucl}$ is the Euclidean metric
has been studied  by the second author and F. Cuccu 
in \cite{culoi}. There they proved the following.

\vskip 0.3cm

\noindent
{\bf Theorem A} 
 {\em Let   $g$
be a  $g_{eucl}$-balanced metric (of height one)
on ${\complex}^n$. If 
$\Phi$  is rotation invariant
then (up to holomorphic isometries)
$g=g_{eucl}$.}

\vskip 0.3cm

In this paper we are concerned with the   $g_{eucl}$-balanced metrics $g$ on
the unit disk ${\cal D}=\{z\in {\complex}\ |\  |z|^2<1\}$, where 
$g_{eucl} =dz\otimes d\bar z$
is  the standard Euclidean metric
on ${\complex}$.
In this case, 
the Hilbert space
${\cal H}$
consists of all holomorphic functions
$\varphi:{\cal D}\rightarrow {\complex}$ such that
$$\int_{{\cal D}}e^{-\Phi}|\varphi|^2\frac{i}{2}dz\wedge d\bar z<\infty ,$$
where $\Phi$ is a \K\ potential for $g$.
Therefore ${\cal H}$ is the weighted Bergman space $L_{hol}^2({\cal D}, e^{-\Phi}\omega_{eucl})$ on ${\cal D}$
with weight $e^{-\Phi}$.
Notice that when  $g=g_{hyp}=\frac{dz\otimes d\bar z}{(1-|z|^2)^2}$ is the hyperbolic metric
on ${\cal D}$,  then  $\Phi (z)=-\log (1- |z|^2)$
is a \K\ potential for $g_{hyp}$ and 
the Hilbert space ${\cal H}=L_{hol}^2({\cal D}, e^{-\Phi}\omega_{eucl})$  consists
of  holomorphic functions  $f$ on ${\cal D}$ such that
$\int_{{\cal D}}(1-|z|^2) |f|^2\frac{i}{2}dz\wedge d\bar z<\infty $. It is easily seen 
 that $\sqrt{\frac{(j+1)(j+2)}{\pi}}z^{j}, j=0,\dots$
is an orthonormal basis of ${\cal H}$.
 The map (\ref{psiglob}), in this case,   is given
by:
$$\varphi_{\Phi}:{\cal D}\rightarrow
{\complex}P^{\infty}:z\mapsto
[\dots,  \sqrt{\frac{(j+1)(j+2)}{\pi}}z^{j}, ,\dots].$$
Thus,
$$\varphi_{\Phi}^*g_{FS}=\frac{i}{2}
\partial\bar\partial\log [\frac{1}{\pi}\sum_{j=0}^{+\infty}
(j+1)(j+2)|z|^{2j} ]=
\frac{i}{2}\partial\bar\partial\log\frac{1}{(1-|z|^2)^3}=3\omega_{hyp}$$
and so $g_{hyp}$ is a $g_{eucl}$-balanced  (even balanced) metric
of height $\alpha =3$.
Notice  that the function $\Phi=-\log (1-|z|^2)$ is a radial function  and $h(x)=e^{-\Phi(z)}=1-|z|^2, x=|z|^2, $ is an entire analytic function   defined on all 
${\real}$.

\vskip 0.3cm

The following theorem, which is the main result of this paper, shows that 
the hyperbolic metric on the unit disk can be characterized by the previous data.

\begin{teor}\label{mainteor}
Let $g$ be a \K\  metric on the unit disk  ${\cal D}$.
Assume that $g$ admits a 
 (globally) defined
\K\ potential $\Phi$ which is radial
and 
such that the function   $h(x)= e^{-\Phi (z)}$,  $x=|z|^2$,  extends to  a (real valued)  entire analytic function
on ${\real}$.
If the metric $g$ is  $g_{eucl}$-balanced of height $3$, 
 then 
$g=g_{hyp}$.
\end{teor}

The  proof of Theorem \ref{mainteor} is based on
the following    characterization of the function $f(x)=1-x$
very interesting on its own sake.

\begin{lemma}\label{mainprop}
  Let $\lambda$ be a positive real number, and let $f \colon \mathbb R \to
  \mathbb R$ be an entire analytic function such that $f(x) > 0$ for all $x
  \in (0,1)$. Define
  \begin{equation}\label{condition0}
  I_j = \int_0^1 f(t) \, t^j \, dt
  \quad
  \mbox{for $j \in \mathbb N$}.
\end{equation}
If the series $\displaystyle \sum_{j = 0}^{+\infty} \frac{\, x^j \,}{\, I_j
  \,}$ converges for every $x \in (-1,1)$, and if
\begin{equation}\label{condition}
\frac{\, 2 \, \lambda^2 \,}{\, f^3(x) \,}
=
\sum_{j = 0}^{+\infty} \frac{\, x^j \,}{\, I_j \,}
\qquad
\mbox{for all $x \in (0,1)$},
\end{equation}
then $f(x) =\lambda  (1 - x)$ for all $x \in \mathbb R$.
\end{lemma}

\vskip 0.3cm

Despite  the very natural statement
the proof of Lermma \ref{mainprop} is far to be trivial, 
being  based on a careful analysis 
of the behaviour
of  $f(x)$ and its derivatives  as $x\rightarrow 1^-$.

\vskip 0.3cm

In view of this lemma
the authors believe the validity 
of the following conjecture 
which could be an important   step
towards the classification  
of $g_{eucl}$-balanced metrics of height $\alpha$ on the complex hyperbolic space
${\complex}H^n$, namely the unit ball  $B^n\subset {\complex}^n$
equipped with the hyperbolic form $\omega_{hyp}=-\frac{i}{2}\partial\bar\partial\log (1-\|z\|^2)$,
$z\in B^n$.

\vskip 0.3cm

\noindent
{\bf Conjecture:} 

\noindent
{\em Fix  a positive integer $n$ and let 
$$D_n=\{x=(x_1, \dots , x_n)\in {\real}^n|\ 0<x_1+\cdots +x_n<1, \ x_j>0\}.$$
Suppose that there exists an integer  $\alpha> n+1$, 
a positive real number $\lambda$ and  an entire analytic function
$f:{\real}^n\rightarrow {\real}$ such that $f(x) > 0$ for all $x\in D_n$ 
such that 
$$\frac{(\alpha -1)\cdots (\alpha -n)\lambda^2}{f^{\alpha}(x)}=
\sum_{J}\frac{x^{J}}{I_{J}(\alpha)},\   \forall J=(j_1, \dots , j_n)\in {\natur}^n,$$
where 
$$I_{J}(\alpha)=\int_{D_n}f^{\alpha -(n+1)}(x)
x^Jdx_1\cdots dx_n.$$
Then $f(x)=\lambda(1-x_1-\cdots -x_n)$.}

\vskip 0.3cm

Notice that Lemma \ref{mainprop} shows the validity of the  previous conjecture
for $n=1$ and  $\alpha =3$.

\begin{remar}\rm
The studies of  balanced metrics  on the unit ball $B^n\subset {\complex}^n$
is far  more complicated that one of studying the 
$g_{eucl}$-balanced metrics 
(we refer the reader to a recent paper of Miroslav Engli\v{s}
\cite{meweig} for the study of radial balanced metrics on $B^n$).
The situation  is similar in the compact case 
where  there are no obstructions for the existence of $g_0$-balanced metrics
(where $g_0$ is  a fixed metric)
on a given integral  \K\ class  of  a compact  complex  manifold $M$  while the existence of balanced metric on $M$ is subordinated  to the existence of  a constant scalar curvature metric in that class (cf.    \cite{arlfirst} and \cite{bly}).
\end{remar}

\begin{remar}\rm
Lemma \ref{mainprop}
should be compared with  the following   characterization of the exponential
function due to Miles and Williamson \cite{mw}  which is the main
tool  in \cite{culoi} in order   to prove  
Theorem A: {\em  
let 
$f(x)=\sum_j b_jx^{j}$
be an entire function on ${\real}$ such that 
$b_0=1,\  b_j>0, \ \forall j\in {\natur},$
and
$$ \int_{{\real}}\frac{b_j t^{j}}{f(t)}\,dt=1,\ \forall j\in {\natur}, $$
then $f(x)=e^x$.}
\end{remar}

\vskip 0.5cm

The paper contains 
another  section where we  prove
Lemma \ref{mainprop} and  Theorem  \ref{mainteor}.

\section{Proof of  the main results }\label{bal}
In the proof of  Lemma \ref{mainprop} we need the
the following elementary result.
\begin{lemma}\label{O}
Let $r_0 \in \mathbb N$. If a sequence $\{ c_j \}$ satisfies $c_j =
O(j^{r_0})$ as $j \to +\infty$, then the power series $\sum_{j =
0}^{+\infty} \, c_j \, x^j$ converges in the interval $(-1, 1)$
to a function $S(x)$ such that
$S(x) = O((1-x)^{-r_0 - 1})$ as $x \to 1^-$.
\end{lemma}
\dimostr
If $r_0 = 0$ then the conclusion follows from the
definition of the symbol~$O$ and the fact that $\sum_{j = 0}^{+\infty} \,
x^j = (1 - x)^{-1}$. If, instead, $r_0 > 0$ then the conclusion follows
similarly after the observation that $a_j = O((j + 1) \cdot \dots \cdot (j +
r_0))$ and $\sum_{j = 0}^{+\infty} \, (j + 1) \cdot \dots \cdot (j + r_0) \,
x^j = r_0! \, (1 - x)^{-r_0 - 1}$.
\fdim

\subsection{Proof of Lemma \ref{mainprop}}

\noindent
By replacing $f(x)$ with $\lambda \, f(x)$ we may assume
$\lambda = 1$. 
Unless otherwise stated, the variable~$x$ ranges in the
interval~$(0,1)$. The starting idea of the proof of Lemma
\ref{mainprop}  is the following. From
the Taylor series of $f(x)$ at $x_0 = 1$,
\begin{equation}\label{Taylor}
f(x)
=
\sum_{k = 0}^{+\infty}
\frac{\, (-1)^k \, f^{(k)}(1) \,}{k!} \, (1 - x)^k
,
\end{equation}
we obtain an asymptotic estimate of the left-hand side of~(\ref{condition})  (with $\lambda =1$)
as $x \to 1^-$. Moreover, by repeatedly integrating by parts we obtain, for
every $j, k_0 \in \mathbb N$
\begin{equation}\label{parts}
I_j
=
\textstyle
\sum\limits_{k = 0}^{k_0}
\frac{(-1)^k \, f^{(k)}(1)}
{\, (j + 1) \cdot \dots \cdot (j + k + 1) \,}
+
\frac{(-1)^{k_0 + 1}}{\, (j + 1) \cdot \dots \cdot (j + k_0 + 1) \,}
\,
\int_0^1 f^{(k_0 + 1)}(t) \, t^{j + k_0 + 1} \, dt
.
\end{equation}
Passing to the reciprocal $1/I_j$ and using Lemma~\ref{O} we obtain an
asymptotic estimate of the right-hand side of~(\ref{condition})  (with $\lambda =1$). Since
equality holds, we subsequently determine $f(1), f'(1),f''(1)$.
Then, the proof is concluded by means of a more sophisticated argument.

\goodbreak\noindent Step 1:   $f(1) = 0$ and $f'(1) = -1$. Denote by $k_0 \in \mathbb
N$ the smallest natural number such that $f^{(k_0)}(1) \ne 0$. By
(\ref{Taylor}) we get $f(x) = \frac1{\, k_0 ! \,} \, (-1)^{k_0} \,
f^{(k_0)}(1) \, (1 - x)^{k_0} \, (1 + O(1-x))$. In the sequel we will make
often use of the following elementary expansion:
\begin{equation}\label{elementary}
(1 + t)^{p} = 1 + p t + O(t^2)
\mbox{ as $t \to 0$}, \ p\in {\real}
,
\end{equation}
which implies, in particular, $(1 + O(1 - x))^{-3} = 1 + O(1 - x)$. Taking
this into account, we deduce
\begin{equation}\label{estimate.1}
\frac{\, 2  \,}{\, f^3(x) \,}
=
\frac{2 \, (k_0!)^3 \, (1 + O(1-x))}
{(-1)^{k_0} \, (f^{(k_0)}(1))^3 \, (1 - x)^{3 k_0}}
\end{equation}
Since we are assuming $f^{(k)}(1) = 0$ for $k < k_0$, and since the integral
in~(\ref{parts}) tends to zero at least as fast as $1/j$ as $j \to +\infty$,
we may write
$$
I_j
=
\frac{(-1)^{k_0} \, f^{(k_0)}(1)}
{\, (j + 1) \cdot \dots \cdot (j + k_0 + 1)}
\,
(1 + O(1/j))
,
$$%
which in turn,  by (\ref{elementary}),  implies
$$
\frac1{\, I_j \,}
=
\frac{\, (j + 1) \cdot \dots \cdot (j + k_0 + 1)}
{(-1)^{k_0} \, f^{(k_0)}(1)}
+ O(j^{k_0})
.
$$%
Taking Lemma~\ref{O} into account, multiplication by $x^j$ followed by
summation over $j$ yields
$$
\sum_{j = 0}^{+\infty} \frac{\, x^j \,}{\, I_j \,}
=
\frac{(k_0 + 1)!}
{\, (-1)^{k_0} \, f^{(k_0)}(1) \, (1 - x)^{k_0 + 2} \,}
+
O((1 - x)^{-k_0 - 1})
\hbox{ as $x \to 1^-$}
.
$$%
By comparing the last equality with~(\ref{estimate.1}) it follows that $k_0$
must satisfy $3 k_0 = k_0 + 2$, and therefore $k_0 = 1$. This implies $f(1)
= 0$ and  $(f'(1))^3 =
f'(1)$. Since $f(x) > 0$ for $x \in (0,1)$, $f'(1)$~must be negative and we
conclude $f'(1) = -1$.

\goodbreak\noindent Step 2: $f''(1) = 0$. By Taylor expansion we have $f(x) = \, (1 -
x) \, [\, 1 + \frac1{\, 2 \,} \, f''(1) \, (1 - x) + O((1 - x)^2) \,]$.
Using~(\ref{elementary}) we get
\begin{equation}\label{expansion.2}
\frac{\, 2 \,}{\, f^3(x) \,}
=
\frac2{\, (1 - x)^3 \,}
-
\frac{\, 3 \, f''(1) \,}{\, (1 - x)^2 \,}
+
O((1 - x)^{-1})
.
\end{equation}
Choosing $k_0 = 2$ in~(\ref{parts}) and arguing as before, we also find
$1/I_j = (j + 1)(j + 2) - (j + 1) \, f''(1) + O(1)$ and therefore
by Lemma~\ref{O}
$$
\sum_{j = 0}^{+\infty} \frac{\, x^j \,}{\, I_j \,}
=
\frac2{\, (1 - x)^3 \,}
-
\frac{f''(1)}{\, (1 - x)^2 \,}
+
O((1 - x)^{-1})
.
$$%
By comparing the last estimate with~(\ref{expansion.2}) for $x \to 1^-$ we
deduce $f''(1) = 0$.

\vskip 0.3cm
At this point one could try to obtain the higher order  derivatives  $f^{(k)}(1)$, $k\geq 3$, 
as in Steps 1 and 2.
Unfortunately this does not work.  Indeed  one can easily verify  that 
by iterating the previous procedure
one gets  $f^{(k)}(1)$, $k\geq 4$ in terms of $f^{(3)}(1)$
but the latter remains undetermined. 
In order to overcome this problem  
notice that  the previous steps imply that 
the function
\begin{equation}\label{z}
z(x):=\frac{\, 2 \,}{\, f^3(x) \,}
-\frac2{\, (1 - x)^3 \,}
-
\frac{\, f'''(1) \, }{1 - x}
\end{equation}
is real analytic  
in a neighbourhood of $x=1$.
 Indeed,  we have $f(x) = (1 - x) \, [1 - \frac16
\, f'''(1) \, (1 - x)^2 + (1 - x)^3 \, \varphi(x)]$ for an entire analytic
function~$\varphi(x)$. Furthermore,  $(1 + t)^{-3} = 1 - 3t + t^2 \,
\psi(t)$,  where $\psi(t)$ is analytic for $t \in (-1, +\infty)$
and the claim follows.

Further,   by~(\ref{condition})  (with $\lambda=1$), 
$z(x)$ admits the following expansion  around the origin
$z(x)=\sum_{j = 0}^{+\infty} \, a_j \, x^j$, 
where
\begin{equation}\label{aj}
a_j = 1/I_j - (j + 1) (j + 2) - f'''(1),
\ \ \mbox{for} \ j
\in \mathbb N.
\end{equation}
The proof  of the lemma will be completed by showing that $z(x)$
vanishes identically.
Indeed this is equivalent to 
\begin{equation}\label{a_j=0}
1/I_j  = (j + 1) (j + 2) + f'''(1),\ \ 
\mbox{for} \ j
\in \mathbb N, 
\end{equation}
which plugged  into~(\ref{condition}) (with $\lambda =1$) gives
\begin{equation}\label{last_representation}
\frac2{\, f^3(x) \,}
=
\frac2{\, (1 - x)^3 \,}
+
\frac{\, f'''(1) \,}{\, 1 - x \,}.
\end{equation}
This  shows that $f(1 - t)$ is an odd function of $t$ and therefore
$f^{(4)}(1) = 0$. Taking this into account, and using~(\ref{parts}) with
$k_0 = 4$ we obtain
$$
I_j
=
\frac1{\, (j + 1) \, (j + 2) \,}
-
\frac{f'''(1)}{\, (j + 1) \cdot \dots \cdot (j + 4) \,}
+
O(j^{-6})
,
$$%
which in turn implies $1/I_j = (j + 1) \, (j + 2) + f'''(1) - 4 \, f'''(1) /
j + O(j^{-2}) $. By comparing the last expansion with~(\ref{a_j=0}) we
deduce $f'''(1) = 0$. This and~(\ref{last_representation}) imply $f(x) = 1-
x$ and this concludes the proof of the lemma.

\vskip 0.3cm
In order to prove that the   
the sequence $\{a_j\}$ vanishes identically we need the following steps.
\vskip 0.3cm

Step 3. For every integer $k_1$ there exists a rational function
$Q_{k_1}(j)$ such that
\begin{equation}\label{a_n}
a_j
=
Q_{k_1}(j) + O(j^{-k_1})
\mbox{ as $j \to +\infty$}.
\end{equation}
Observe, firstly, that if~(\ref{a_n}) holds for a particular $k_1 =
\overline k_1$, then it also holds for every $k_1 < \overline k_1$ with
$Q_{k_1} = Q_{\overline k_1}$. Hence, it suffices to prove~(\ref{a_n}) for
$k_1 \ge 1$. Letting $k_0 = k_1 + 2$ in~(7) we obtain: 
$$\textstyle
I_j
=
[ (j + 1) \, (j + 2)]^{-1}
[1 +
\tilde Q_{k_1}(j) 
+
O(j^{-k_1 - 2})],$$
where
$$\textstyle
\tilde Q_{k_1}(j)=
\sum\limits_{k = 3}^{k_1 + 2}
\frac{(-1)^k \, f^{(k)}(1)}{\, (j + 3) \cdot \dots \cdot (j + k + 1)}
=\frac{- f'''(1)}{(j+3)(j+4)}+O(j^{-3}).$$
Therefore
$$1/I_j 
=
\textstyle
(j + 1) \, (j + 2)
\,
[
1
+
\hat Q_{k_1}(j)
+
O(j^{-k_1 - 2})],$$
where 
$\hat Q_{k_1}(j)=-\frac{\tilde Q_{k_1}(j)}{1+\tilde Q_{k_1}(j)}.$

 Letting $Q_{k_1}(j)
= (j + 1) \, (j + 2) \, \hat Q_{k_1}(j) - f'''(1)$, the claim follows by
the definition~(\ref{aj}) of~$a_j$.
In the next step we will need the observation that 
\begin{equation}\label{Qk1}
Q_{k_1}(j)=O(j^{-1}).
\end{equation}

Step 4. The sequence $\{ a_j \}$ defined before tends to zero faster than
every rational function of $j$, namely $a_j = O(j^{-k_1})$ as $j \to
+\infty$ for every integer~$k_1$. This is proved by showing that for every
$k_1 \ge 2$ and every rational function~$Q_{k_1}$ satisfying~(\ref{a_n}) we
have $Q_{k_1}(j) = O(j^{-k_1})$. Suppose that this is not the case. Then, by (\ref{Qk1}),
there exist positive integers $d < k_1$ and a rational function~$Q_{k_1}$
satisfying~(\ref{a_n}) such that the limit $\lim_{j \to +\infty} j^d \,
Q_{k_1}(j)$ is a finite~$c \ne 0$. This and~(\ref{a_n}) imply $a_j = c \,
j^{-d} + O(j^{-d - 1})$. Now recall that the sum of the series $\sum_{j =
  1}^{+\infty} \, x^j/j$ is the unbounded function $-\log(1 - x)$, while the
series $\sum_{j = 1}^{+\infty} \, x^j/j^2$ converges to a bounded function
in the interval $[-1, 1]$.
By comparison with these elementary series it follows that the \hbox{$(d -
  1)$-th} derivative of $\sum_{j = 0}^{+\infty} \, a_j \, x^j$ is unbounded
for $x$ close to~$1^-$. But this is impossible because the last series
converges to $z(x)$, which is analytic in a neighbourhood of $x = 1$. This
contradiction shows that $Q_{k_1}(j) = O(j^{-k_1})$ and the claim follows.

\goodbreak Step 5. The sequence $\{a_j\}$ is identically zero. Define $w_j =
[(j + 1) (j + 2) + f'''(1)] \, I_j - 1$. Since $w_j = - I_j \, a_j$ and
$I_j$ is positive, it suffices to show that $w_j = 0$ for all~$j \in \mathbb
N$. This is achieved by representing $w_j$ as the limit $\lim_{k_0 \to
  +\infty} S_{k_0}(j)$ of the sum $S_{k_0}(j)$ defined below, an then by
showing that $S_{k_0}(j)$ is infinitesimal as $k_0 \to +\infty$. Taking into
account that $\int_0^1 (1 - t)^k \, t^j \, dt = k! \, j! / (j + k + 1)!$,
multiplication of~(\ref{Taylor}) by $t^j$ followed by termwise integration
over the interval $(0,1)$ yields
$$
I_j
=
\textstyle
\sum\limits_{k = 0}^{+\infty}
\frac{(-1)^k \, f^{(k)}(1)}
{\, (j + 1) \cdot \dots \cdot (j + k + 1) \,}
$$%
Notice that it makes sense to integrate over the interval $(0,1)$ since, by
assumption, $f$ is entire (cf.~Remark~\ref{entire}). Since $f(1) = f''(1)
= 0$ and $f'(1) = -1$, the preceding formula leads to
$$
w_j
=
\textstyle
\frac{f'''(1)}{\, (j + 1) (j + 2) \,}
+
\sum\limits_{k = 3}^{+\infty}
\frac{(-1)^k \, [(j + 1) (j + 2) + f'''(1)] \, f^{(k)}(1)}
{\, (j + 1) \cdot \dots \cdot (j + k + 1) \,}
.
$$%
For $k_0 \ge 3$ we may write $w_j = S_{k_0}(j) + R_{k_0}(j)$, where the
partial sum $S_{k_0}(j)$ and the remainder $R_{k_0}(j)$ are given by
\begin{eqnarray}\label{def.S}
S_{k_0}(j)
&=&
\textstyle
\frac{f'''(1)}{\, (j + 1) (j + 2) \,}
+
\sum\limits_{k = 3}^{k_0}
\frac{(-1)^k \, [(j + 1) (j + 2) + f'''(1)] \, f^{(k)}(1)}
{\, (j + 1) \cdot \dots \cdot (j + k + 1) \,}
,
\\ \noalign{\medskip}\nonumber
R_{k_0}(j)
&=&
\textstyle
\sum\limits_{k = k_0 + 1}^{+\infty}
\frac{(-1)^k \, [(j + 1) (j + 2) + f'''(1)] \, f^{(k)}(1)}
{\, (j + 1) \cdot \dots \cdot (j + k + 1) \,}
.
\end{eqnarray}
By~(\ref{parts}), the remainder $R_{k_0}(j)$ also admits the following
representation:
$$
R_{k_0}(j)
=
\textstyle
\frac{\, (-1)^{k_0 + 1} \, [(j + 1) (j + 2) + f'''(1)] \,}
{(j + 1) \cdot \dots \cdot (j + k_0 + 1)}
\int_0^1 f^{(k_0 + 1)}(t) \, t^{j + k_0 + 1} \, dt
,
$$%
which shows that $R_{k_0}(j) = O(j^{-k_0})$ as $k_0 \to +\infty$.
Furthermore, since $w_j = - I_j \, a_j$ and $I_j$ is bounded, by  Step~4 
we have, in particular,
$w_j = O(j^{-k_0})$. It follows that $S_{k_0}(j) = O(j^{-k_0})$ and
by~(\ref{def.S}) we may write
\begin{equation}\label{S}
S_{k_0}(j)
=
\frac{P^1_{k_0}(j)}{\, (j + 1) \cdot \dots \cdot (j + k_0 + 1) \,}
,
\end{equation}
where $P^1_{k_0}(j) = m_{k_0} \, j + q_{k_0}$ is a convenient polynomial of
degree $\deg P^1_{k_0} \le 1$ in the variable $j$. In order to show that
$S_{k_0}(j)$ is infinitesimal as $k_0 \to +\infty$ we have to investigate
the coefficients $m_{k_0}, q_{k_0}$. 
Observe, firstly, that from~(\ref{def.S}) we get
$$
\textstyle
S_{k_0 + 1}(j)
=
S_{k_0}(j)
+
\frac{(-1)^{k_0 + 1} \, [(j + 1) (j + 2) + f'''(1)] \, f^{(k_0 + 1)}(1)}
{\, (j + 1) \cdot \dots \cdot (j + k_0 + 2) \,}
.
$$
This and~(\ref{S}) yield
$$
\textstyle
\frac{P^1_{k_0 + 1}(j)}{\, (j + 1) \cdot \dots \cdot (j + k_0 + 2) \,}
=
\frac{P^1_{k_0}(j)}{\, (j + 1) \cdot \dots \cdot (j + k_0 + 1) \,}
+
\frac{(-1)^{k_0 + 1} \, [(j + 1) (j + 2) + f'''(1)] \, f^{(k_0 + 1)}(1)}
{\, (j + 1) \cdot \dots \cdot (j + k_0 + 2) \,}
.
$$%
By summation of the two rational functions in the right-hand side of the
last equality, and since the coefficients of $j^2, j, j^0$ in the numerator
must equal the corresponding ones in the left-hand side, we deduce
\begin{eqnarray*}
0\quad &=& m_{k_0} + (-1)^{k_0 + 1} \, f^{(k_0 + 1)}(1)
,
\\
m_{k_0 + 1} &=& (k_0 - 1) \, m_{k_0} + q_{k_0}
,
\\
q_{k_0 + 1} &=& (k_0 + 2) \, q_{k_0} - [2 + f'''(1)] \, m_{k_0}
.
\end{eqnarray*}
Since the series~(\ref{Taylor}) converges together with all its derivatives
at $x = 0$, it follows that for every~$h \in \mathbb N$ we have $m_{k_0} = o
(k_0! \, k_0^{-h})$ as $k_0 \to +\infty$. The same holds for $q_{k_0}$
because $q_{k_0 + 1} = (k_0 + 2) \, [m_{k_0 + 1} - (k_0 - 1) \, m_{k_0}] -
[2 + f'''(1)] \, m_{k_0}$. Hence, by~(\ref{S}), it follows that $S_{k_0}(j)
\to 0$ as $k_0 \to +\infty$ and therefore $w_j = 0$ for all $j$, as claimed.
\fdim

\begin{remar}\rm\label{entire}
The assumption in Lemma \ref{mainprop} that $f$ is an
entire analytic function  can be relaxed. If equality~(\ref{Taylor}) holds in
the interval $(-\varepsilon, \, 2 + \varepsilon)$ for some $\varepsilon >
0$, then the same proof shows that $f(x) = \lambda \, (1 - x)$ in that
interval. 
\end{remar}


\subsection{ Proof of Theorem \ref{mainteor}}\label{secteor}
Since the function   $h(x)=e^{-\Phi (z)},\  x=|z|^2$,  extends to all real numbers it follows 
 that  $e^{-\Phi (z)}$ does not blow up at the boundary of ${\cal D}$. This  implies  
 that   the monomials $z^j$, $j=0, 1\dots$
 are an orthogonal basis  of ${\cal H}=L_{hol}^2({\cal D}, e^{-\Phi}\omega_{eucl})$.
 Hence the sequence $\sqrt{b_j}z^j, j=0,\dots$, with
$$b_j=\left(\int_{\cal D}e^{-\Phi}|z|^{2j}\frac{i}{2}dz\wedge d\bar z\right) ^{-1},$$ 
 is an orthonormal basis of  ${\cal H}$ and 
 the \K\  metric $g$
is $g_{eucl}$-balanced
of height $3$ iff
$$\frac{i}{2}\partial\bar\partial\log (\sum_{j=0}^{\infty}b_j|z|^{2j})=3\omega=3\frac{i}{2}\partial\bar\partial\Phi.$$
This implies that  the function 
$ \Phi(z)-\log(\sum_{j=0}^{\infty}b_j|z|^{2j})^{\frac{1}{3}}$
is a radial harmonic function on ${\cal D}$
and hence equals  a constant, say $\Phi_0$.
By setting  $f(x)=(\sum_{j=0}^{\infty}b_{j}x^{j})^{-\frac{1}{3}}, \ x\in [0, 1)$, 
and by the definition of the  $b_j$'s
one then  gets
\begin{equation}\label{alm}
e^{-\Phi_0}b_{j}\int_{{\cal D}}f(|z|^2)|z|^{2{j}}
\frac{i}{2}dz\wedge d\bar z
=1,\ \forall
j\in {\natur} .
\end{equation}
Observe  that  again the assumption that $h(x)=e^{-\Phi (z)}, \ x=|z|^2$,
extends to  an entire analytic function on ${\real}$
implies the same property for 
$f(x)=h(x)e^{\Phi_0}$. 
By passing to  polar coordinates
$z=\rho e^{i\theta}, \ 
\rho\in [0, +\infty),\ \theta\in [0, 2\pi)$
and by the change of variables $t=\rho^2$
one obtains:
$$\pi e^{-\Phi_0}b_{j}\int_{0}^1
f(t)t^{{j}}dt
=1,\ \forall
j\in {\natur}.$$
By setting
 $\lambda^2=\frac{\pi e^{-\Phi_0}}{2}$, 
 $I_j=\frac{1}{2\lambda^2b_j}$
and by the definition of $f(x)$
one gets (\ref{condition0})
and (\ref{condition}).
Therefore,   by Lemma \ref{mainprop},  $f(x)=\lambda (1-x)$,
i.e.
$\Phi (z)=\Phi_0-\log\lambda-\log (1-|z|^2)$, 
which implies $\omega =\frac{i}{2}\partial\bar\partial\Phi =
-\frac{i}{2}\partial\bar\partial\log (1-|z|^2)=\omega_{hyp}$
 and this concludes the proof of the theorem.
\fdim

\small{}

\end{document}